\documentclass[11pt,reqno,oneside]{article}
\usepackage{amsmath}
\usepackage{amssymb}
\usepackage[dvips]{graphicx}

\usepackage{pstricks}

\newtheorem{prop}{Proposition}
\newtheorem{coro}{Corollary}

\title{Unbalanced subtrees in binary rooted ordered and un-ordered trees}
\author{Filippo Disanto\thanks{fdisanto@uni-koeln.de}}
\date{}

\begin{document}

\maketitle

\abstract{Binary rooted trees, both in the ordered and in the un-ordered case, are well studied structures in the field of combinatorics. The aim of this work is to study particular patterns in these classes of trees. We consider completely unbalanced subtrees, where unbalancing is measured according to the so-called Colless's index. The size of the biggest unbalanced subtree becomes then a new parameter with respect to which we find several enumerations.}

\section{Introduction}

The aim of this work is to study particular unbalanced patterns in rooted binary trees, both in the ordered and un-ordered case. More precisely, we consider a new statistic on trees. We are interested in the size of the biggest subtree having the \emph{caterpillar} property.

Caterpillars have already been considered in the case of \emph{coalescent} trees, see for example the interesting work of Rosenberg~\cite{rosenberg}. In particular, in a genetic population framework, when trees are used to represent ancestry relations among individuals, the presence of a caterpillar subtree often indicates phenomena such as  \emph{natural selection}.
 
Subtrees structures have already been considered \cite{lara} and \cite{rowland}. However, to our knowledge, enumerative properties of caterpillar subtrees have never been investigated as done in this work. 

In Section~\ref{ordered} we start by giving some basic definitions. We then provide the enumeration of ordered rooted binary trees of a given size having the biggest caterpillar subtree of size less than, greater than or equal to a fixed integer $k$. Furthermore, we provide the expected value of the size of the biggest caterpillar subtree when ordered trees of size $n$ are uniformly distributed and $n$ is large.

In Section~\ref{permut} we see how caterpillar subtrees correspond to patterns extracted from $132$-avoiding permutations. The resulting characterization is interesting and will represent a starting point for further studies on sub-structures of permutations. 

Finally, in Section~\ref{unordered} we study caterpillars realized in un-ordered binary rooted trees. The resulting approach is similar to the one used in the ordered case and it provides asymptotic formulas for the probability of a tree of a given size with "small" caterpillar subtrees.

\section{Caterpillars in ordered rooted binary trees}\label{ordered}

\subsection{Definitions}\label{def}

Ordered rooted binary trees are enumerated with respect to the size, i.e., number of leaves, by the well known sequence of \emph{catalan} numbers corresponding to entry $A000108$ in \cite{sloane}. The generating function of catalan numbers is denoted by $C(x)$ and it looks as follows

\begin{equation}\label{catralan}
C(x)=\frac{1- \sqrt{1-4x}}{2}.
\end{equation}

We denote the class of ordered rooted binary trees by $\mathcal{T}$ while $\mathcal{T}_n$ denotes the subset of $\mathcal{T}$ made of those elements having size $n$. In Section~\ref{ordered} we use the term tree referring to ordered binary rooted trees.

\begin{figure}
\begin{center}
\includegraphics*[scale=.4,trim=0 0 0 0]{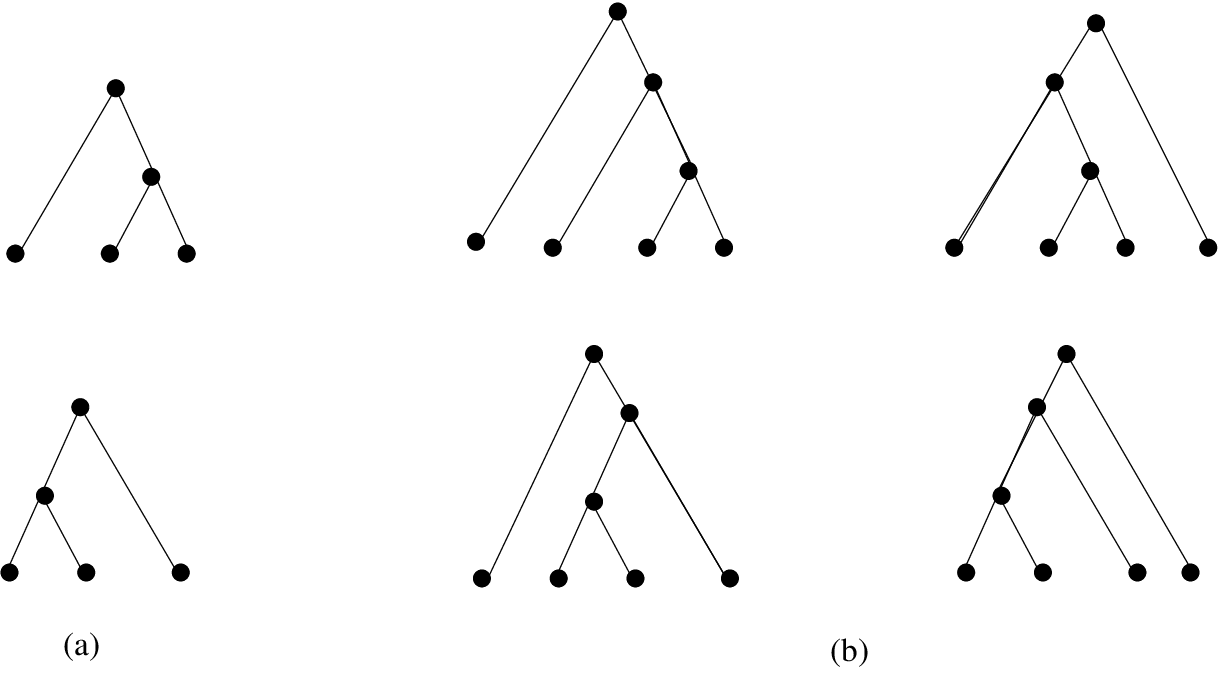}
\end{center}
\caption{(a)  caterpillars of size 3; (b) caterpillars of size 4.}\label{cat}
\end{figure}

We define a  tree in $\mathcal{T}_n$ to be a \emph{caterpillar} of size $n$ if each node is a leaf or it has  at least one leaf as a direct descendant. See for example Fig.~\ref{cat}~(a)~(b).

In addition, caterpillars can be characterized by the fact that they are the most unbalanced trees. As a measure of tree imbalance we take the following index.
Given a tree $t$ and a node $i$, let $t_l(i)$ (resp. $t_r(i)$) be the left (resp. right) subtree of $t$ determined by $i$. We define $$\Delta_t(i) = |\text{size}(t_l(i)) - \text{size}(t_r(i))|.$$ 

If $t \in \mathcal{T}_n$, its \emph{Colless}'s index (see \cite{slatkin}) is defined as

$$\frac{2}{(n-2)(n-1)} \cdot \sum_{i \, \text{node of} \, t} \Delta_t(i).$$ 

The Colless's index is considered as a measure of tree imbalance. Its value ranges between $0$ and $1$, where $0$ corresponds to a completely balanced tree while $1$ to an unbalanced one. 

Based on the previous definitions, a tree of size $n>2$ is a caterpillar if and only if its Colless's index is $1$.

\begin{figure}
\begin{center}
\includegraphics*[scale=.4,trim=0 0 0 0]{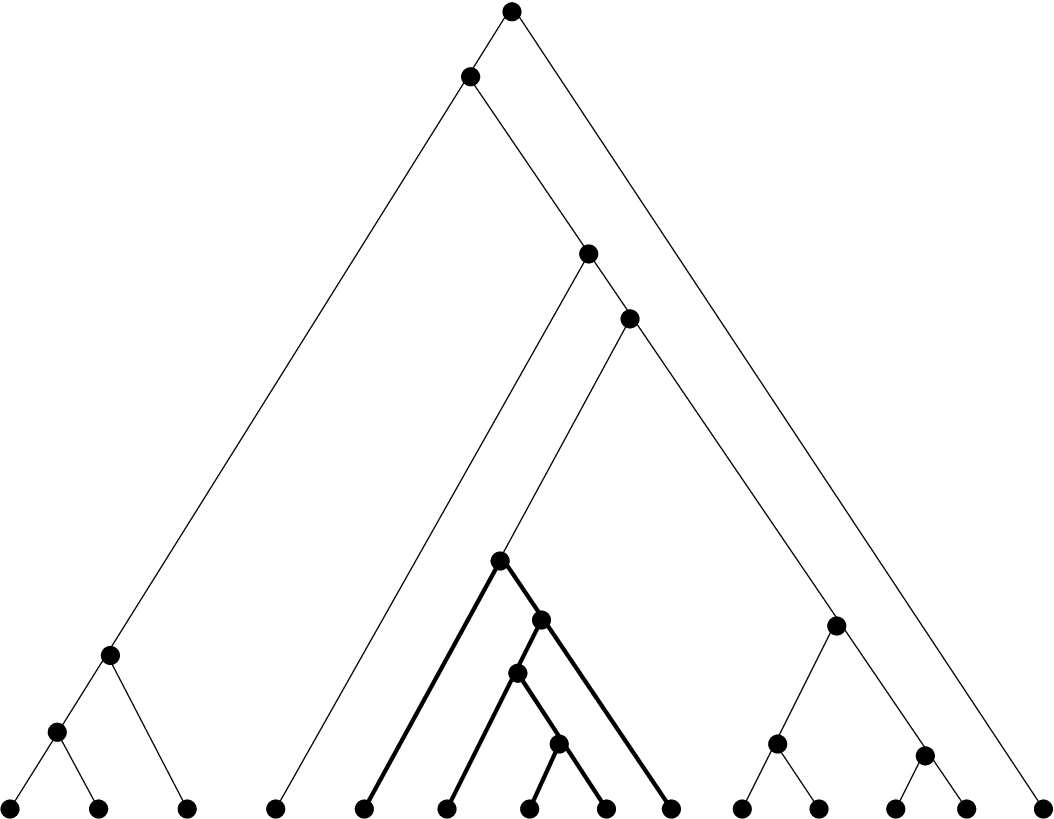}
\end{center}
\caption{A tree having $\gamma$ parameter equal to 5. The biggest caterpillar is highlighted.}\label{dim}
\end{figure}

If $t \in \mathcal{T}_n$, we define $\gamma(t)$ as the size of the biggest caterpillar 
which can be seen as a subtree of $t$. We observe that, if $n>1$, then $\gamma(t) \geq 2$. In Fig.~\ref{dim} we show a tree having $\gamma=5$.

\subsection{A recursive construction for the size of the biggest caterpillar subtree}\label{gf}

Let $F^{-}_k(x)$ be the ordinary generating function which gives the number of trees having $\gamma$ parameter \emph{at most} equal to $k \geq 2$.

It is easy to see that $F^{-}_k$ satisfies the equation 

\begin{equation} \label{rec}
F^{-}_k = x + (F^{-}_k)^2 - 2^{k-1} x^{k+1}.
\end{equation}

Indeed a tree $t$ having $\gamma(t) \leq k$ has either size one or it is built appending two trees $t_1$ and $t_2$, with $\gamma(t_1) \leq k$ and $\gamma(t_2) \leq k$, to the root of $t$. We must  exclude the case in which one among $t_1$ and $t_2$ has size $1$ and the other one is a caterpillar of size $k$. Since there are exactly $2^{k-2}$ caterpillars of size $k$ the previous formula follows.  

From (\ref{rec}) we obtain 

\begin{equation}\label{ordinato}
F^{-}_k(x)=\frac{1-\sqrt{1 - 4x + 2^{k+1}x^{k+1}}}{2}.
\end{equation}

Then, considering $F^{+}_k = C(x) - F^{-}_{k-1}(x)$, one has the number of trees having $\gamma \geq k$ while, taking $F_k = F^{-}_k(x) - F^{-}_{k-1}(x)$, one can compute the number of trees having $\gamma = k$. The following table shows the first coefficients of the Taylor expansion of $F^{-}_k$, $F^{+}_k$ and $F_k$ when $k=5$. 

\begin{center} 
\begin{tabular}{|c|cccccccccc|}
\hline
k=5 & 1 & 2 & 3 & 4 & 5 & 6 & 7 & 8 & 9 & 10 \\\hline
$F^{-}_k$ & 1 & 1 & 2 & 5 & 14 & 26 & 100 & 333 & 1110 & 3742 \\ 
$F^{+}_k$ & 0 & 0 & 0 & 0 & 8 & 16 & 48 & 160 & 560 & 1952 \\ 
$F_k$ & 0 & 0 & 0 & 0 & 8 & 0 & 16 & 64 & 240 & 832 \\\hline
\end{tabular} 
\end{center}

Note that the sixth coefficient of $F_k$ is $0$. Indeed, as the reader can easily check, there is no tree of size $k+1$ having $\gamma$ parameter equal to~$k$.

We conclude this section observing that none of the sequences corresponding to $F^{-}_k$, $F^{+}_k$ and $F_k$ is present at the moment in \cite{sloane} with other combinatorial interpretations.

\subsubsection{Asymptotic growth}

The function $F_k(x)$ has its singularities for the solutions of the equation  
$1 -4x +2^{k+1}x^{k+1} = 0$. By \emph{Pringsheim}'s theorem (see \cite{ancomb}) we can assume, for our purposes, that the dominant singularity of $F_k(x)$ corresponds to the positive real solution of $1 -4x +2^{k+1}x^{k+1} = 0$ which is closer to the origin. 
Let $\rho_k$ be this solution. We observe that, when $k$ increases, $\rho_k$ approaches $1/4$. 
In order to prove this claim we remark that, for $k \geq 2$, we have 

\begin{equation}\label{est2}
\frac{1}{4} < \rho_k < \frac{2}{5}.
\end{equation}

Indeed, this can be shown by considering the polynomial $$y=1 -4x +2^{k+1}x^{k+1}$$ which satisfies $y(1/4)>0$, $y(2/5)<0$ and it is such that $y \geq 1-4x > 0$ for $0 \leq x < 1/4$. 
We now proceed by \emph{bootstrapping} (see \cite{ancomb}). Writing the defining equation for $\rho_k$ as $$x=\frac{1}{4}(1+2^{k+1}x^{k+1})$$ and making use of (\ref{est2}) yields $$\frac{1}{4}\left( 1+\frac{1}{2^{k+1}}\right) < \rho_k < \frac{1}{4}\left( 1+\left( \frac{4}{5} \right)^{k+1}\right)$$ which is sufficient to prove that $\rho_k \rightarrow 1/4$.    

A further iteration of the previous inequality shows that 
\begin{small}
$$\rho_k < \frac{1}{4} \left( 1+2^{k+1}\left( \frac{1}{4} \left(1+ \left( \frac{4}{5} \right)^{k+1} \right) \right)^{k+1} \right)$$
\end{small}
which, considering that $$\left(1+ \left( \frac{4}{5} \right)^{k+1} \right)^{k+1} = 1+(k+1)\left(  \frac{4}{5} \right)^{k+1} + \mathcal{O}\left( k^2\left( \frac{4}{5}  \right)^{2k}  \right), $$ gives 
\begin{small}
$$ \rho_k < \frac{1}{4} \left( 1+ \frac{1}{2^{k+1}} + (k+1)\left(\frac{2}{5} \right)^{k+1}      +\mathcal{O}\left(k^2 \left(\frac{8}{25}\right)^k\right)\right).$$
\end{small}

Thus $$\rho_k - \frac{1}{4} - \frac{1}{2^{k+3}} <  \frac{1}{4}(k+1)\left( \frac{2}{5}\right)^{k+1} + \mathcal{O}\left(k^2 \left(\frac{8}{25}\right)^k\right),$$ which means $$\rho_k = \frac{1}{4} + \frac{1}{2^{k+3}} + \mathcal{O}\left( k  \left( \frac{2}{5}\right)^{k}  \right).$$

In the following table we show the first approximated values of $\rho_k$.

\begin{center} 
\begin{tabular}{|c|c|}
\hline
$\rho_2$ & 0.3090169 \\ 
$\rho_3$ & 0.2718445 \\ 
$\rho_4$ & 0.2593950 \\
$\rho_5$ & 0.2543301 \\ 
$\rho_6$ & 0.2520691 \\ 
$\rho_7$ & 0.2510085 \\\hline
\end{tabular} 
\end{center}

Now observe that, for a given constant $a$, we can always write
$$1 -4x +2^{k+1}x^{k+1} = (a -x)(4 - 2^{k+1} \sum_{i=0}^{k} a^{i} x^{k-i}) + 1 -4 a +2^{k+1} a ^{k+1}.$$
Substituting $a$ with $\rho_k$ one has
$$1 -4x +2^{k+1}x^{k+1} = (\rho_k -x)(4 - 2^{k+1} \sum_{i=0}^{k} \rho_k^{i} x^{k-i}).$$
If we now set $$B(x) = 4 - 2^{k+1} \sum_{i=0}^{k} \rho_k^{i} x^{k-i},$$
by standard asymptotic calculations (see \cite{ancomb}) we obtain for large $n$
\begin{align} \label{asi}
[x^{n}] F^{-}_k &\sim  \frac{1}{4} \sqrt{ \frac{B(\rho_k) \rho_k}{\pi n^3}} \left( \frac{1}{\rho_k} \right) ^n \\\nonumber
&= \frac{1}{4} \sqrt{\frac{4 \rho_k -(k+1)2^{k+1} \rho_k^{k+1}}{\pi n^3}} \left( \frac{1}{\rho_k} \right)^{n}. \nonumber
\end{align}

We can apply formula (\ref{asi}) to provide the asymptotic behaviour of trees with no caterpillar of size $3$. Caterpillars with three leaves are also called \emph{pitchforks} in \cite{rosenberg}.

\begin{prop}
The number of pitchfork-free trees of size $n$ is given by $[x^n]F^{-}_2$ and it satisfies asymptotically the following relation:
$$[x^n] F^{-}_2(x) \sim \frac{1}{4} \sqrt{\frac{4 R -24 R^{3}}{\pi n^3}} \left( \frac{1}{R} \right)^{n},$$
where $R=\frac{1}{4}(\sqrt{5}-1)=0.3090169$.
\end{prop}

When $n=100$ the ratio between $[x^{100}] F^{-}_2$ and its approximation is $0.9933$.

\subsection{The average size of the biggest caterpillar subtree}

In this section we determine the value $E_n(\gamma)$ which denotes the average of $\gamma(t)$ when $t \in \mathcal{T}_n$.

As shown in Section~\ref{gf}, when $k>0$,  $F^{-}_k(x)$ gives the number of trees having $\gamma$ at most $k$. Indeed, also in the case $k=1$, we have $F^{-}_1=(1-\sqrt{1-4x+4x^2})/2=x$ which represents the unique caterpillar of size $1$. 

Furthermore consider $f^{(n)}_k=[x^n]F^{-}_k(x)$ and analogously let us denote by $C^{(n)}=[x^n]C(x)$  the $n$-th catalan number. We can write the desired average value as 

\begin{align}\nonumber
E_n(\gamma)&=\frac{1 f^{(n)}_1 + \sum_{k \geq 1} (k+1)(f^{(n)}_{k+1} - f^{(n)}_{k})}{C^{(n)}} \\\nonumber
&= \frac{-f^{(n)}_1-...-f^{(n)}_{n-1}+nf^{(n)}_n+ \sum_{k \geq n} (k+1)(f^{(n)}_{k+1} - f^{(n)}_{k})}{C^{(n)}} \\\nonumber
&= \frac{-f^{(n)}_1-...-f^{(n)}_{n-1}+nC^{(n)}+ \sum_{k \geq n} (C^{(n)}-f^{(n)}_k)}{C^{(n)}} \\\nonumber
&= \frac{\sum^{n-1}_{k=1} (C^{(n)}-f^{(n)}_k) +C^{(n)}+\sum_{k \geq n} (C^{(n)}-f^{(n)}_k)}{C^{(n)}} \\\nonumber
&= \frac{C^{(n)}+ \sum_{k \geq 1} (C^{(n)}-f^{(n)}_k)}{C^{(n)}} \\\nonumber
&= 1 + \frac{\sum_{k \geq 1} (C^{(n)}-f^{(n)}_k)}{C^{(n)}}. \nonumber
\end{align}

In the previous calculation we rely on the fact that $f^{(n)}_k=C^{(n)}$ for $k \geq n$.
We now focus our attention on the generating function $U(x)$ which is defined as
$$U(x)=\sum_{k \geq 1} (C(x)-F^{-}_k(x))= \frac{1}{2} \sum_{k \geq 1} \left( \sqrt{1-4x+2^{k+1} x^{k+1}} -\sqrt{1-4x} \right).$$
Near the dominant singularity $x=1/4$ we can substitute the term $\sqrt{1-4x+2^{k+1} x^{k+1}}$ with the corresponding $\sqrt{1-4x+\frac{1}{2^{k+1}}}$. The effect of the substitution in the sum can be measured, expanding up to the first order, as
\begin{small}
$$\left|\sum_{k\geq 1}\sqrt{1-4x+2^{k+1} x^{k+1}} - \sum_{k\geq 1}\sqrt{1-4x+\frac{1}{2^{k+1}}}\right|\sim \left| x-\frac{1}{4} \right|\sum_{k \geq 1} (k+1)\left(\frac{1}{2^{k-1}}\right)^{1/2},$$
\end{small}
where $\sum_{k \geq 1} (k+1)\left(\frac{1}{2^{k-1}}\right)^{1/2}\simeq 15.071$. By the mentioned substitution we obtain the generating function $V(x)$ whose coefficients grow asymptotically like those of $U(x)$,
$$V(x)=\frac{1}{2} \sum_{k \geq 1} \left(\sqrt{1-4x+\frac{1}{2^{k+1}}} -\sqrt{1-4x} \right).$$
Considering the $n$-th coefficient of $V(x)$ we define the sequence $(g_n)_{n\geq 1}$ as
$$g_n=\frac{[x^n]V(x)}{C^{(n)}}=\sum_{k\geq 1} \left( 1-\left( 1+\frac{1}{2^{k+1}} \right)^{-n+1/2} \right)$$
which gives an asymptotic approximation of $E_n(\gamma)-1$. 

By a \emph{Poisson/Mellin-transform} approach (see \cite{spanoski} for details) we can now further investigate the growth of  the coefficients $g_n$ when $n$ is large. By a Poisson-transform we reduce the problem to the asymptotic analysis of a \emph{harmonic} sum which is then studied using Mellin transforms. 

Setting 
$$C=\sum_{k\geq 2}  \left( 1-\left( 1+\frac{1}{2^{k}} \right)^{1/2} \right)\simeq -0.24056,$$
we write
$$g_n=C+\sum_{k\geq 2} \left( 1+\frac{1}{2^k} \right)^{1/2}\left( 1-\left( 1+\frac{1}{2^k} \right)^{-n} \right)= C + h_n.$$ 
If $H(x)$ is the exponential generating function of the sequence $(h_n)_{n\geq 1}$, we compute the associated Poisson-transform $\tilde{H}(x)$  yielding
$$\tilde{H}(x)=\sum_{n\geq 0} h_n \frac{x^n}{n!}\exp(-x)=\sum_{k\geq 2}(1+2^{-k})^{1/2}\left( 1-\exp \left( \frac{-x}{2^k+1} \right) \right).$$ We are now interested in the behaviour of $\tilde{H}(x)$ when $x \rightarrow \infty$. Indeed, for $n$ large, $h_n$ is approximated by $\tilde{H}(n)$.  
Observe that $\tilde{H}(x)$ is a harmonic sum, i.e., it is of the form $$\tilde{H}(x)=\sum_k \lambda_k \tilde{h}(\mu_k \cdot x).$$ In the open strip of complex numbers $s=\alpha +i \beta$ such that $-1<\alpha<0$, the associated Mellin-transform is
\begin{align}\nonumber
\mathcal{M}(\tilde{H};s)=& \sum_{k\geq 2} \frac{\lambda_k}{\mu_k^s} \cdot\int_0^{\infty}\tilde{h}(x)x^{s-1}dx= -\Gamma(s)\sum_{k\geq 2} (1+2^{-k})^{1/2}(1+2^k)^s \\\nonumber
=& -\Gamma(s)\sum_{k\geq 2} 2^{ks}[ ( 1+2^{-k} )^{s+1/2}-1 ] -\Gamma(s)\sum_{k\geq 2} (2^s)^k \\\label{transf}
=&  -\Gamma(s)\sum_{k\geq 2} 2^{ks}[ ( 1+2^{-k} )^{s+1/2}-1 ] +\Gamma(s)\frac{4^s}{2^s-1}. \\\nonumber
\end{align}

The behaviour of $\tilde{H}(x)$ for $x$ large can be obtained by the analysis of singular expansions of (\ref{transf}) to the right of the strip $-1<\Re(s)<0$. The transform can be analytically continued in $0<\Re(s)\leq M$ for any $M>0$, then the poles of interest are just those at $s=0$ (double pole) and at $s=\chi_k=\frac{2k\pi \, i}{\log 2}$ for $k\in \mathbb{Z}\setminus \{0\}$. 

The singular expansion of $\Gamma(s)$ at  $s=0$ looks like $\Gamma(s)\sim \frac{1}{s}~ -~\eta$ (with $\eta \simeq 0.57721$ being Euler's constant) and similarly one has $\frac{1}{2^s-1}\sim \frac{1}{s\log 2}-\frac{1}{2}$. Furthermore we need to consider the expansion $4^s\sim 1+s\log 4$. Putting all together we obtain the expansion of (\ref{transf}) near the double pole $s=0$ as
$$\frac{1}{s^2\log 2}+\frac{1}{s}\left( -\frac{\eta}{\log 2}+\frac{3}{2}+C \right)+\mathcal{O}(1).$$ 

Instead, near $s=\chi_k$, we expand (\ref{transf}) as $$\frac{\Gamma(\chi_k)}{\log 2\cdot (s-\chi_k)}+\mathcal{O}(1),$$ given that $4^{\chi_k}=1$.

\begin{figure}
\begin{center}
\includegraphics*[scale=.8,trim=0 0 0 0]{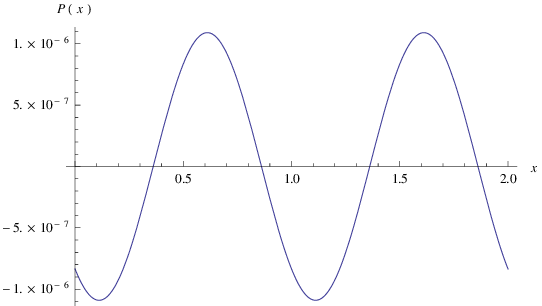}
\end{center}
\caption{Fluctuations determined by $P(x)$.}\label{fluttu}
\end{figure}

Each pole $s_0$ contributes to the asymptotic of $\tilde{H}(x)$ with a term determined by the following rule (see again \cite{spanoski}):
$$\frac{d}{(s-s_0)^{k+1}} \rightarrow -\frac{(-1)^k d}{k!} \cdot x^{-s_0}(\log x)^k.$$
In this way, when $x\rightarrow \infty$, we find for any $M>0$
$$\tilde{H}(x) = \log_2(x)+\frac{\eta}{\log 2}-\frac{3}{2}-C - \frac{P(\log_2 x)}{\log 2} +\mathcal{O}(x^{-M}),$$
where 
\begin{equation}\label{fluttua}
P(x)=\sum_{k\neq 0}\Gamma(\chi_k)\exp(-2k\pi i \cdot x)
\end{equation}
is a function of period $1$ with mean zero and minute fluctuations  bounded by $\max(|P(x)|)\simeq 10^{-6}$, see Fig.~\ref{fluttu}. 

The behaviour of the coefficients $h_n$ can be reached by Depoissonization obtaining
$$h_n=\log_2(n)+\frac{\eta}{\log 2}-\frac{3}{2}-C - \frac{P(\log_2 n)}{\log 2} +\mathcal{O}\left(\frac{\log^*n}{n}\right)$$ from which 
$$E_n(\gamma)=\log_2(n)+\frac{\eta}{\log 2}-\frac{1}{2} - \frac{P(\log_2 n)}{\log 2} +\mathcal{O}\left(\frac{\log^*n}{n}\right).$$

According to the previous calculations we can state

\begin{prop} \label{palla}
When $n$ is large, the expected size of the biggest caterpillar sub-tree in a tree of size $n$ is given by
\begin{equation}\label{befana}
E_n(\gamma)=\log_2(n)+\frac{\eta}{\log 2}-\frac{1}{2} - \frac{P(\log_2 n)}{\log 2} +\mathcal{O}\left(\frac{\log^*n}{n}\right),
\end{equation}
where $\eta \simeq 0.57721$ is Euler's constant and $P(x)$ is a small periodic fluctuation of mean zero defined by (\ref{fluttua}).
\end{prop}

Evaluating the non-fluctuating term we obtain the second row of the following table while in the first row we find, for several values of $n$, the true $E_n(\gamma)$ which is computed by generating functions.

\begin{center} 
\begin{tabular}{|c|ccccc|}
\hline
$n$ & 50 & 100 & 200 & 500 & 1000 \\\hline
$E_n(\gamma)$ & 6.202 & 7.107 & 8.052 & 9.334 & 10.318 \\ 
(\ref{befana}) & 5.976 & 6.976 & 7.976 & 9.298 & 10.298 \\\hline
\end{tabular} 
\end{center}

As a corollary to Proposition~\ref{palla} we can state

\begin{coro}\label{coropalle}
When $n \rightarrow \infty$ we have

$$\frac{\log_2(n)}{E_n(\gamma)} \sim 1.$$

\end{coro}

\section{Caterpillars in permutations $Av(132)$}\label{permut}

In Section~\ref{def} we have introduced  caterpillars as objects related to trees. We know that also the class of permutations avoiding the pattern $132$ is enumerated by catalan numbers. Indeed, one can bijectively map the set $\mathcal{T}_{n+1}$ onto the set $Av_{n}(132)$, where the last symbol denotes the class of permutations of size $n$ avoiding $132$. In particular, we will use a  bijection $\phi: \mathcal{T}_{n+1} \rightarrow Av_{n}(132)$ which works as described below. 
\medskip 

\begin{figure}
\begin{center}
\includegraphics*[scale=.4,trim=0 0 0 0]{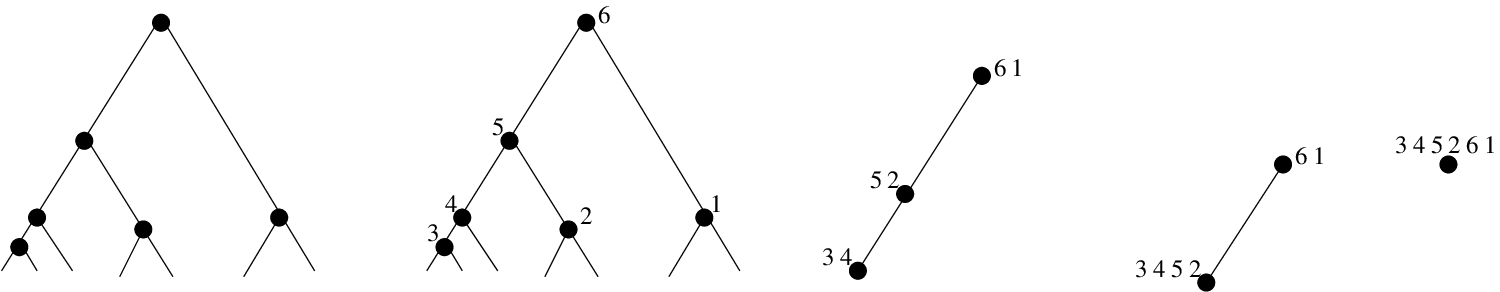}
\end{center}
\caption{The mapping $\phi$.}\label{bij}
\end{figure}

\textbf{Bijection $\phi$.}
Take $t \in \mathcal{T}_{n+1}$ and visit it according to the pre-order traversal. At the same time, starting with the label $n$ for the root, label each node of outdegree two in decreasing order. After this first step one obtains a tree with integers associated with the nodes of outdegree two. Each leaf now collapses to its direct ancestor which takes a new label receiving on the left and right the label of its left and right child respectively. We go on collapsing leaves until we obtain a tree made of one node which is labelled with a permutation of size $n$. See Fig.~\ref{bij} for an instance of this mapping.  
\medskip 

Using $\phi$ one can see how caterpillars are realized in permutations with no $132$ pattern. We need the following definition. Let $\pi=\pi_1 \pi_2 \dots \pi_n$ be a permutation. For a given entry $\pi_i$ we define $r_{\pi}(\pi_i)$  as the set made of those entries $\pi_k$ such that:
\begin{enumerate}
\item[1)] $\pi_k \leq \pi_i$;
\item[2)] all entries of $\pi$ which are placed between $\pi_k$ and $\pi_i$ are less than or equal to $\pi_i$.
\end{enumerate}

\begin{figure}
\begin{center}
\includegraphics*[scale=.4,trim=0 0 0 0]{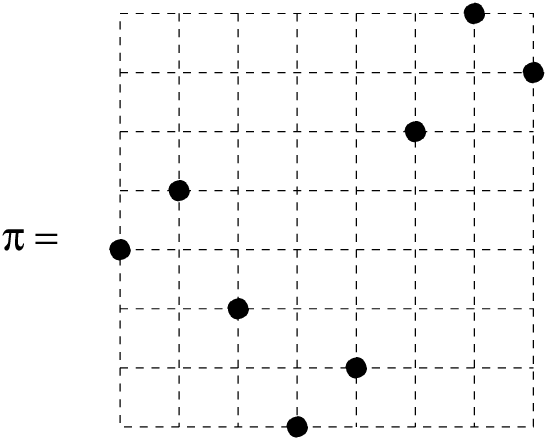}
\end{center}
\caption{The permutation $\pi =  (4 5 3 1 2 6 8 7)$.}\label{permutaz}
\end{figure}

Given $\pi=\pi_1 \pi_2 \dots \pi_n$,  let $\tilde{r}_{\pi}(\pi_i)$ be the permutation one obtains extracting from $\pi$ the elements of $r_{\pi}(\pi_i)$ respecting the order. The set of permutations $\{\tilde{r}_{\pi}(\pi_i)\}_{i=1\dots n}$ is denoted by $\tilde{r}_{\pi}$. As an example, consider the permutation $\pi$ which is shown in Fig.~\ref{permutaz}. In this case $\tilde{r}_{\pi}$ is made of 

\begin{align}\nonumber
\tilde{r}_{\pi}(4) & =  (1), \\\nonumber  
\tilde{r}_{\pi}(5) & =  (4 5 3 1 2), \\\nonumber 
\tilde{r}_{\pi}(3) & =  (3 1 2), \\\nonumber 
\tilde{r}_{\pi}(1) & =  (1), \\\nonumber
\tilde{r}_{\pi}(2) & =  (1 2), \\\nonumber
\tilde{r}_{\pi}(6) & =  (4 5 3 1 2 6), \\\nonumber
\tilde{r}_{\pi}(8) & =  (4 5 3 1 2 6 8 7), \\\nonumber
\tilde{r}_{\pi}(7) & =  (1). \nonumber
\end{align}

The next proposition describes how caterpillars look like in permutations avoiding the pattern $132$. It is interesting to see that the presence of such particular subtrees is linked to the property of avoiding the pattern $231$. 

\begin{prop}\label{chepalle}
If $t \in \mathcal{T}_{n+1}$ and $\phi(t)=\pi=\pi_1 \pi_2 \dots \pi_n$, then the following holds:
\begin{itemize}
\item[i)] caterpillar subtrees of $t$ correspond through $\phi$ to those permutations in $\tilde{r}_{\pi}$ avoiding the pattern $231$;
\item[ii)] $\gamma(t)-1$ corresponds to the size of the biggest permutation in 
\begin{equation} \nonumber
Av(231) \cap \tilde{r}_{\pi}.
\end{equation}
\end{itemize}
\end{prop}
\emph{Proof.}
Label $t$ according to the procedure $\phi$. If a node is labelled with $m$ consider the subtree $t_m$ whose root is $m$. The nodes belonging to $t_m$  form the subsequence of $\pi$ made of the elements of $r_{\pi}(m)$. Therefore, we find the pattern $231$ in $\tilde{r}_{\pi}(m)$ if and only if we can find a node in $t_m$ having two descendants which are not leaves of $t$. It is now sufficient to observe that $t_m$ is a caterpillar if and only if it does not contain such a node.
Summarizing, for every node $m$ of $t$, $t_m$ is a caterpillar subtree of size $k+1$ if and only if $\tilde{r}_{\pi}(m) \in Av_k(231)$.
$\Box$

\bigskip

Using the results of Proposition~\ref{chepalle} as well as those contained in previous sections, we can describe some properties of the permutations in $\tilde{r}_{\pi}$ when $\pi$ avoids the pattern $132.$ Indeed we have the next two corollaries.

\begin{coro}\label{duepalle}
The generating function of the number of permutations $\pi \in Av(132)$ such that all elements in $\tilde{r}_{\pi}$ of size greater than one contain the pattern $231$  is given by $$\frac{F^{-}_2(x)}{x}-1=\frac{1-2x-\sqrt{1 - 4x + 8x^{3}}}{2x}.$$ The first terms of the sequence are: $$1,0,1,2,6,16,45,126,358,1024,2954,
8580,25084,73760,218045.$$   
\end{coro}

{\bf Remark.} Given $\pi=\pi_1 \pi_2 \dots \pi_n$, we say that $\pi_i$ is a \emph{valley} when $\pi_{i-1}$ and $\pi_{i+1}$ (if they exist) are greater than $\pi_i$; while $\pi_i$ is said to be a \emph{peak} if both $\pi_{i-1}$ and $\pi_{i+1}$ exist and  $\pi_{i-1} < \pi > \pi_{i+1}$. In this sense, those permutations $\pi$ considered in Corollary~\ref{duepalle} can be characterized, among those in $Av(132)$, by the fact that each entry $\pi_i$  either is a valley or has $\tilde{r}_{\pi}(\pi_i)$ containing at least one peak. 
We also observe that sequence A025266 of \cite{sloane} provides the same list of numbers given by the previous corollary. The mentioned sequence also enumerates Motzkin paths with additional constraints.

\medskip

Finally we state the following result which can be deduced from Corollary~\ref{coropalle}.

\begin{coro}
If $\pi \in Av(132)$ has size $n$, the expected size of the biggest permutation in $Av(231) \cap \tilde{r}_{\pi}$ is asymptotic to $\log_2(n).$ 
\end{coro}

\section{Caterpillars in un-ordered rooted binary trees}\label{unordered}

In the previous sections we have focused our attention on the presence of caterpillar subtrees in ordered rooted binary trees. As a second step we would like to investigate the un-ordered case. Let us start recalling some basic enumerative properties.

Un-ordered rooted binary trees are enumerated with respect to the size, i.e., number of leaves, by the sequence $w_1,w_2,...,w_n,...$ of the so-called \emph{Wedderburn-Etherington} numbers. This sequence corresponds to entry $A001190$ of \cite{sloane}. The corresponding generating function  $W(x)$ is defined implicitly by the following equation

\begin{equation}\nonumber
W(x)=x+ \frac{1}{2}W(x)^2 +\frac{1}{2}W(x^2).
\end{equation}

The asymptotic behaviour of $(w_n)_{n>0}$ is given by

\begin{equation}\label{wedderburn}
w_n \sim \frac{\lambda}{2\sqrt{\pi}} n^{-3/2}\rho^{-n},
\end{equation}

where $$\frac{\lambda}{2\sqrt{\pi}}=0.3187766259 \mathrm{\, \, and \,\,} \rho=0.40269750367.$$ See for example \cite{ancomb}.
 
\bigskip
 
The class of un-ordered rooted binary trees is denoted by $\mathcal{W}$ while $\mathcal{W}_n$ represents the subset of $\mathcal{W}$ whose elements have size $n$. In the present section we use the word tree to refer to un-ordered binary rooted trees.

The definition of a \emph{caterpillar} tree in $\mathcal{W}$ is the same as in the ordered case. We have to pay attention to the fact that now, due the un-ordering constraint, the number of different caterpillar trees of fixed size is one (see again Fig.~\ref{cat}~(a)~(b)). Analogously to Section~\ref{ordered} we define the parameter $\gamma(t)$ as the size of the biggest caterpillar subtree of the tree $t$.

Let $W_k(x)$ be the ordinary generating function which gives the number of trees having $\gamma$ \emph{at most} equal to $k>0$.

One can see that, similarly to the functions $F^{-}_k$ of Section~\ref{ordered}, $W_k$ satisfies the equation 

\begin{equation} \label{ritorno}
W_k(x) = x + \frac{1}{2}W_k(x)^2 + \frac{1}{2}W_k(x^2) - x^{k+1}.
\end{equation}

The generating function $W_k(x)$ has radius of convergence $\rho_k$ which is at least $\rho=0.402\dots$ and (for $k\geq 2$) at most $1/2$. The latter bound can be observed from a comparison with the solution to the functional equation $g=x+x^2/2+g^2/2-x^3$. Indeed the solution $g$ counts those trees with no caterpillar of size $3$ where each tree of size greater than two such that its left root sub-tree is isomorphic to the right root sub-tree is counted $1/2$. From (\ref{ritorno}) we obtain 

\begin{equation} \label{piove}
W_k(x)= 1 - \sqrt{1 - 2\left( x + \frac{W_k(x^2)}{2} - x^{k+1} \right)} = 1 - \sqrt{1 - 2 \varphi_k(x)}
\end{equation}

and we see that $\rho_k$ corresponds to the smallest positive solution of $\varphi_k(x)=1/2$. 

The function $\varphi_k(x)$ is analytic in the disc $|x|<\rho_k^{1/2}$ which then contains the one determined by $\rho_k$. Expanding $\varphi_k(x)$ near $x=\rho_k$ gives

$$\varphi_k(x) = \varphi_k(\rho_k) + \left(1-\rho_k^{k}-k\rho_k^{k}+\rho_k W^{'}_k(\rho_k^2) \right) (x-\rho_k) + \mathcal{O}((x-\rho_k)^2)$$

and plugging into (\ref{piove}) we obtain the singular expansion of $W_k(x)$ at $x=\rho_k$ as

\begin{small}
\begin{align}\nonumber
W_k(x) & =  1 - \sqrt{2\rho_k-2\rho_k^{k+1}-2k\rho_k^{k+1}+2\rho_k^2 W^{'}_k(\rho_k^2) } \cdot \sqrt{1-\frac{x}{\rho_k}} + \mathcal{O}((x-\rho_k)^{3/2}) \\\label{pasqui}
&= 1 - \lambda_k \cdot \sqrt{1-\frac{x}{\rho_k}}+ \mathcal{O}((x-\rho_k)^{3/2}). 
\end{align}
\end{small}

Starting from (\ref{pasqui}) and perfoming a standard singularity analysis we obtain the number of trees of size $n$ having gamma parameter at most $k$. This number is denoted by $w_{n,k}$ and asymptotically we have

\begin{equation}\label{asinosardo}
w_{n,k} = \frac{\lambda_k}{2\sqrt{\pi}} n^{-3/2} \rho_k^{-n} + \mathcal{O}(n^{-5/2}\rho_k^{-n}).
\end{equation}  

We now proceed as described in \cite{ancomb} providing a procedure which numerically approximates the constants $\rho_k$ and $\lambda_k$ which are involved in (\ref{asinosardo}). The accuracy of the approximations depends on a parameter $m$ which is here taken as $m=10$.

Once we have fixed $k$, we can compute the numbers $w_{1,k}, w_{2,k},..., w_{m,k}$ by recursively applying (\ref{ritorno}). The values for $k=1...5$ are listed in the table below.

\begin{center} 
\begin{tabular}{|c|cccccccccc|}
\hline
$k$ & $w_{1,k}$ & $w_{2,k}$ & $w_{3,k}$ & $w_{4,k}$ & $w_{5,k}$ & $w_{6,k}$ & $w_{7,k}$ & $w_{8,k}$ & $w_{9,k}$ & $w_{10,k}$ \\\hline
1 & 1 & 0 & 0 & 0 & 0 & 0 & 0 & 0 & 0 & 0 \\ 
2 & 1 & 1 & 0 & 1 & 1 & 2 & 3 & 6 & 10 & 19 \\ 
3 & 1 & 1 & 1 & 1 & 2 & 4 & 7 & 14 & 27 & 55 \\
4 & 1 & 1 & 1 & 2 & 2 & 5 & 9 & 19 & 37 & 78 \\ 
5 & 1 & 1 & 1 & 2 & 3 & 5 & 10 & 21 & 42 & 89 \\\hline
\end{tabular} 
\end{center}

Using these entries we define 
$$\tilde{\varphi_k}(x) = x - x^{k+1} + \frac{1}{2}\sum_{i=1}^{m} w_{i,k} x^{2i}$$
and we consider $\rho_k$ approximated by the smallest positive solution $\tilde{\rho_k}$ of $\tilde{\varphi_k}(x)=1/2$ . We can estimate also $W^{'}_k(\rho_k^2)$ as $ \sum_{i=1}^{m} i \cdot w_{i,k} \cdot \tilde{\rho_k}^{2i-2}$ finding an approximation $\tilde{\lambda_k}$ for $\lambda_k$. We observe that increasing the precision $m>10$ - as we will see in the next paragraph - does not change the first five (resp. four) digits of $\tilde{\rho_k}$ (resp. $\tilde{\lambda_k}$). It is then reasonable to assume that, up to five (resp. four) digits,  $\tilde{\rho_k}=\rho_k$ (resp. $\tilde{\lambda_k}=\lambda_k$).  

The results for $k=2...5$ are listed in the following table. In the last column we find the ratio beween the real value of $w_{50,k}$ and the one given by (\ref{asinosardo}) calculated using $\tilde{\rho_k}$ and $\tilde{\lambda_k}$.

\begin{center} 
\begin{tabular}{|c|ccc|}
\hline
$k$ & $\tilde{\rho_k} $ & $\tilde{\lambda_k}/(2\sqrt{\pi})$ & $w_{50,k}/$(\ref{asinosardo}) \\\hline 
2 & 0.46745 & 0.2789 & 1.008\\ 
3 & 0.42291 & 0.2991  & 1.009\\
4 & 0.41001 & 0.3089 & 1.010\\ 
5 & 0.40550 & 0.3139 & 1.011\\\hline
\end{tabular} 
\end{center}

Formula (\ref{asinosardo}), together with (\ref{wedderburn}), gives the probability of a tree of size $n\gg1$ having no caterpillar of size greater than $k$. As an example take $n=100$ and $k=5$, then $$\frac{w_{100,5}}{w_{100}} \sim \frac{0.3139}{0.3187} \times \left( \frac{0.40550}{0.40269} \right)^{-100}=0.984 \times (1.006)^{-100} \sim 0.5.$$ 

Roughly speaking 50\% of trees of size $100$ have no caterpillar of size greater than $5$.  

\paragraph{Small caterpillars.} In order to use the asymptotic result (\ref{asinosardo}), we list, in this final section, the values of $\rho_k$ and $\lambda_k/(2\sqrt{\pi})$ for $2\leq k \leq 10$. Furthermore, we do this more precisely than before. Indeed we choose $m=30$ for our approximations which gives, at least, $10$ digits of accuracy with respect to the exact values. In this way we will be able to find the probabilities $w_{n,k}/w_n$ when $n$ is large and $k$ is small (i.e. less than or equal to $10$). The table below shows the values we are interested in. 

\bigskip

\begin{center} 
\begin{tabular}{|c|cc|}
\hline
$k $ & $\rho_k$      & $\lambda_k/(2\sqrt{\pi})$   \\\hline 
2    & 0.4674554078  & 0.2789408958\\ 
3    & 0.4229139375  & 0.2991123692  \\
4    &  0.4100112389 & 0.3089581337\\ 
5    &  0.4055024052 & 0.3139472095  \\
6    & 0.4038017227  & 0.3164492710\\ 
7    & 0.4031375239  & 0.3176775180  \\
8    & 0.4028738458  & 0.3182668950\\ 
9    &  0.4027683607 & 0.3185438777  \\
10   & 0.4027260095  & 0.3186717321 \\\hline
\end{tabular} 
\end{center}

Using these values, we plot in Fig.~6 the probability of a tree with at least one caterpillar of size greater than $k$, i.e. $1 - (w_{n,k}/w_n)$, for large $n$ and $k=3,4,5,8$.

\begin{figure}
\begin{center}\label{ma}
\begin{tabular}{|c|c|}\hline
\includegraphics*[angle=0,scale=.6,trim=0 0 0 0]{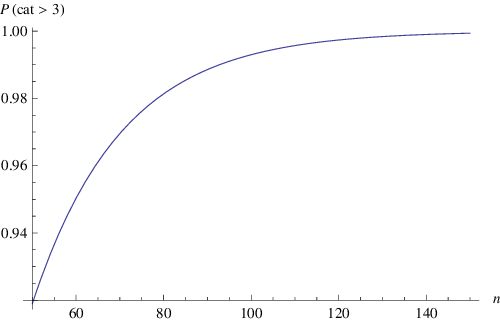} & \includegraphics*[angle=0,scale=.6,trim=0 0 0 0]{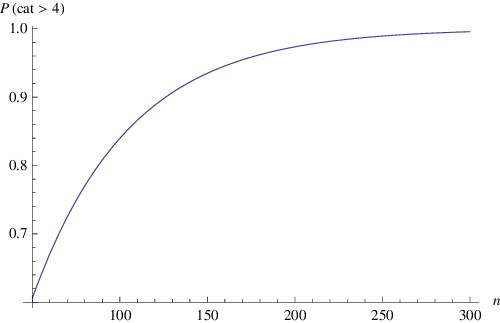} \\ \hline
\includegraphics*[angle=0,scale=.6,trim=0 0 0 0]{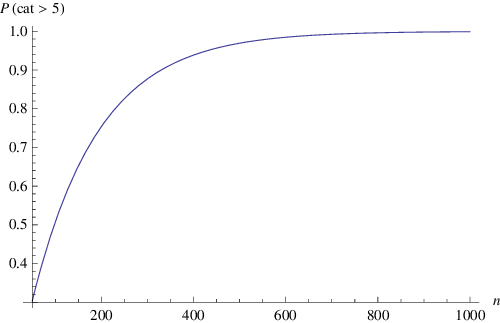} & \includegraphics*[angle=0,scale=.6,trim=0 0 0 0]{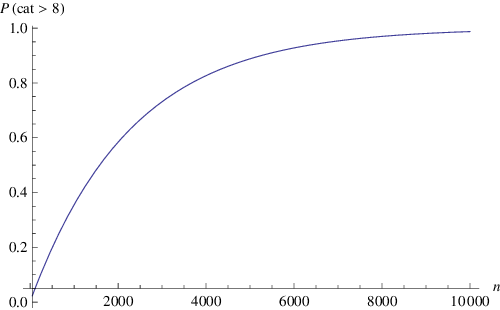} \\ \hline
\end{tabular}
\end{center}
\caption{The probability of a tree of size $n$ with at least one caterpillar of size $k=3,4,5,8$.}
\end{figure}

\section*{Acknowledgements}

We would like to thank an anonymous referee for several suggestions that contributed to improve this manuscript.

\end{document}